\documentstyle[twoside,amssymb]{article}

\makeatletter
\def\cicy#1(#2|#3)#4{\left(\matrix{#2}\right|\!\!
		     \left|\matrix{#3}\right)^{{#4}}_{#1}}

\def\P#1{{\bf P}^{#1}}

\def\C{{\Bbb C}}

\def\N{{\Bbb N}}
\def\Q{{\Bbb Q}}
\def\R{{\Bbb R}}
\def\P{{\Bbb P}}
\def\Z{{\Bbb Z}}

\newcommand{\vskp}{\vspace{5mm}}

\newtheorem{theo}{Theorem}[section]

\newtheorem{rem}[theo]{Remark}
\newtheorem{coro}[theo]{Corollary}
\newtheorem{conj}[theo]{Conjecture}
\newtheorem{prop}[theo]{Proposition}
\newtheorem{exam}[theo]{Example}
\newtheorem{dfn}[theo]{Definition}
\newtheorem{ques}[theo]{Question}

\begin{document}

\begin{center}
\Large
{\bf STRINGY  HODGE  NUMBERS  OF   
     VARIETIES  WITH \\ 
\vskp
GORENSTEIN CANONICAL  SINGULARITIES
}
\normalsize

\vskp
\vskp
\vskp
{V. V. BATYREV}

\vskp
{\it Mathematisches Institut, Universit\"at T\"ubingen, \\
Auf der Morgenstelle 10, D-72076 T\"ubingen, Germany
\\E-mail: batyrev@bastau.mathematik.uni-tuebingen.de}
\end{center}

\begin{center}
{\bf Abstracts}
\end{center}

\begin{quote}
{\footnotesize 
We introduce the notion of stringy $E$-function for an arbitrary 
normal irreducible algebraic variety $X$ with at worst log-terminal 
singularities. We prove some basic properties 
of stringy $E$-functions and compute them explicitly 
for arbitrary $\Q$-Gorenstein toric varieties. Using stringy $E$-functions,
we propose  a general method  to define  stringy Hodge numbers 
for projective algebraic varieties with at worst   
Gorenstein canonical singularities. This  allows us to formulate  
the  topological mirror duality test for arbitrary 
Calabi-Yau varieties with canonical singularities.  In Appendix 
we explain non-Archimedian integrals over spaces of arcs. We need these  
integrals for  the proof of the main technical statement 
used 
in the definition of stringy Hodge numbers. 
}
\end{quote}

\section{Introduction}
Topological field theories 
associated with Calabi-Yau manifolds  predict the following 
Hodge-theoretic property of mirror manifolds \cite{W}: 
\bigskip

{\em If  two smooth $d$-dimensional Calabi-Yau manifolds $V$ and $V^*$ 
form a mirror pair, then their Hodge numbers must satisfy the relations
\begin{equation}
h^{p,q}(V) = h^{d-p,q}(V^*), \;\; 0 \leq p,  q \leq d. 
\label{d-rel} 
\end{equation}
In particular, one has the following equality for the Euler numbers
\begin{equation}
e(V) = (-1)^d e(V^*). 
\end{equation}
}

The relation (\ref{d-rel}) is known as a {\em topological mirror symmetry 
test} \cite{M}. The physicists' discovery of  mirror symmetry has been  
supported by   a lot of examples of mirror pairs $(V,V^*)$ 
consisting of Calabi-Yau varieties with Gorenstein abelian quotient 
singularities \cite{GP}. Some further examples of mirror pairs  
$(V,V^*)$ consisting of Calabi-Yau varieties with Gorenstein toroidal 
singularities  have been proposed  in \cite{BA,BS,BB0,LB}. We want to stress 
that a mathematical formulation of the {topological  mirror symmetry test} 
for Calabi-Yau varieties with singularities is a priori not clear,       
because  (\ref{d-rel}) does not 
hold for usual Hodge numbers of $V$ and $V^*$ if at least one of 
these Calabi-Yau varieties is not smooth. 

Let us illustrate  this problem using some simplest examples: 

\begin{exam}
{\rm Let $V \subset \P^4$ be a Fermat quintic 3-fold defined by the 
equation $f(z)= \sum_{i =0}^4 z_i^5 =0$. Physicists predict 
that the  mirror of $V$ is isomorphic to the quotient $V^*= V/G$, where 
$G \cong (\Z/5Z)^3$ is maximal diagonal subgroup of $PSL(5,\C)$ which 
acts trivially on the polynomial $f$. It is easy to check that 
\[ 1=   h^{2,1}(V^*) =  h^{2,1}(V^*) = 
 h^{1,1}(V)=  h^{2,2}(V). \]
However one has 
\[ 1 = h^{1,1}(V^*) = h^{2,2}(V^*) \neq  h^{2,1}(V) =  h^{2,1}(V) =101. \]
In order to obtain the required duality (\ref{d-rel}), one needs to 
construct  a smooth Calabi-Yau $3$-fold ${\widehat{V}}^*$ by resolving  
singularities of $V^*$. It is possible to show that the Hodge numbers
$h^{1,1}({\widehat{V}}^*)= h^{2,2}({\widehat{V}}^*)= 101$ do not  depend on 
the choice of such a desingularization (see \cite{B} for more general 
result).} 
\end{exam}

\begin{exam}
{\rm Let $V \subset \P^5$ be a Fermat octic  4-fold defined by the 
equation $f(z)= z_0^4 + z_1^4 + z_2^8 +  z_3^8 +  z_4^8 +  z_5^8
=0$. Then physicists predict that the  mirror of $V$ is isomorphic 
to the quotient $V^*= V/G$, where 
$G$ is maximal diagonal subgroup of $PSL(6,\C)$ which 
acts trivially on the polynomial $f$. By resolving abelian 
quotient singularities on  $V^*$, one obtains a smooth Calabi-Yau 
$4$-fold ${\widehat{V}}^*$ with the Hodge numbers 
\[ h^{1,1}({\widehat{V}}^*) =  h^{3,3}({\widehat{V}}^*) = 433,\;\; 
  h^{2,2}({\widehat{V}}^*) = 1820.\]
On the other hand, it is easy to check that 
\[ h^{3,1}(V) = h^{1,3}(V) = 433,  h^{2,2}(V) = 1816.\]
We observe that $h^{2,2}({\widehat{V}}^*) -  h^{2,2}(V) = 4 \neq 0$. 
The duality  (\ref{d-rel}) fails, because $V$ has $4$ isolated singular 
points: $ z_0^4 + z_1^4 =0$, $z_2 = z_3=z_4 =z_5=0$. A small neigbourhood 
of each such a point is analytically isomorphic to $\C^4/\pm id$.  
Unfortunately, it is impossible 
to resolve these singularities of $V$ if one wants to keep  
the canonical class trivial.} 
\label{exam2} 
\end{exam} 

The last example shows that we need   
some new Hodge numbers, so called {\em string-theoretic Hodge numbers}
(or simply {\em stringy Hodge numbers}) $h^{p,q}_{\rm st}(X)$  for singular 
Calabi-Yau varieties $X$,  so that for each mirror pair $(V, V^*)$ of  
singular  Calabi-Yau varieties we  would have  
\begin{equation}
h^{p,q}_{\rm st}(V) = h^{d-p,q}_{\rm st}(V^*), \;\; 0 \leq p,  q \leq d. 
\label{d-rel2} 
\end{equation}
(In particular, we would have  $h^{2,2}_{\rm st}(V) = 
h^{2,2}(V) + 4 
=1820$ 	in  \ref{exam2}.)  
\bigskip

Let us fix our terminology. In this paper we want to consider the 
following most general definition of Calabi-Yau varieties 
with singularities:

\begin{dfn} 
{\rm  A normal projective algebraic variety $X$ is called a 
{\bf Calabi-Yau variety} if $X$ has at worst Gorenstein canonical 
singularities and the canonical line bundle on $X$ is trivial. } 
\label{d-cy}
\end{dfn} 

\begin{rem} 
{\rm One usually requires  that the cohomologies  of the structure sheaf 
${\cal O}_X$ of a Calabi-Yau variety $X$ satisfy some  additional vanishing 
conditions 
\[ h^i(X, {\cal O}_X) = 0, \;\;\;\; 0 < i < d= dim\,X. \]
However, in the connection with the mirror symmetry, we don't want to restrict 
ourselves to these additional assumptions. }
\end{rem}

So far (see \cite{BD,BB}), our approach to the topological mirror symmetry 
test was restricted  to the following steps: First, one chooses   
a class  of  ``allowed'' Gorenstein  canonical singularities 
on Calabi-Yau varieties. Second, one defines  the  notion of 
{\em stringy Hodge numbers} $h^{p,q}_{\rm st}$   for Calabi-Yau with the 
allowed singularities. Third, one proves  (\ref{d-rel2}) for already known 
examples mirror pairs $(V,V^*)$ of Calabi-Yau varieties with the allowed 
singularities. For instance, a  definition of stringy Hodge numbers for      
projective algebraic varieties with Gorenstein quotient 
singularities or with Gorenstein toroidal singularities 
has been introduced in \cite{BD}, and 
all steps of the above  program has been succesfully accomplished for 
Calabi-Yau complete intersections in toric varieties in  \cite{BB}.

One could naturally ask whether there exists a way to define  
stringy Hodge numbers for projective algebraic varieties with 
{\em arbitrary} Gorenstein canonical singularities. One of purposes in   
this paper is to  show that in general the answer is ``no''. Nevertheless, 
there always exists a way to formulate the 
topological  mirror 
symmetry test for Calabi-Yau varieties with arbitrary Gorenstein 
canonical singularities. Our main idea consists of an association  with 
every  projective variety $X$ with at worst Gorenstein singularities a 
rational function $E_{\rm st}(X; u,v) \in \Z[[u,v]] \cap \Q(u,v)$, which 
we call {\em stringy $E$-function}. 
If $E_{\rm st}(X; u,v)$ is a polynomial, then the coefficients of this 
polynomial allow us to define the stringy Hodge numbers  
$h^{p,q}_{\rm st} (X)$ by the 
formula 
\[ E_{\rm st}(X; u,v) = \sum_{p,q}  (-1)^{p+q} h^{p,q}_{\rm st}(X)u^p v^q. \]
If $E_{\rm st}(X; u,v)$ is not a polynomial, then we can't define  stringy 
Hodge numbers of $X$ and say that they  
{\em do not exist}. However,  even in this bad situation 
we could think about  $E_{\rm st}(X; u,v)$ as if it were a generating 
function for Hodge numbers of a smooth projective manifold, i.e., as 
\[ \sum_{p,q} (-1)^{p+q} h^{p,q}u^p v^q. \]
In particular, we shall see that the function  $E_{\rm st}(X; u,v)$ 
always satisfies the following Poincar\'e duality: 
\[  E_{\rm st}(X; u,v)=  (uv)^d E_{\rm st}(X; u^{-1},v^{-1}). \]
For this reason, one can always formulate the topological mirror symmetry 
test for a mirror pair $(V,V^*)$ of $d$-dimensional Calabi-Yau varieties 
with arbitrary singularities as the equality 
\[   E_{\rm st}(V; u,v)=  (-u)^{d}E_{\rm st}(V^*; u^{-1},v), \]
whenever the stringy Hodge numbers of $V$ and $V^*$ exist or not. 
This is exactly the formula which has been 
proved  for Calabi-Yau complete intersections in toric varieties \cite{BB}. 

The paper is organized as follows. In section 3 we introduce 
stringy $E$-functions and discuss their properties. In section 4
we compute these functions for toric varieties. In section 5 we 
give further examples of $E$-functions and formulate some natural 
open questions. In Appendix we explain the proof of the main 
technical statement \ref{key} which rely upon a non-Archimedian integration
over the space of arcs. This part of the paper is strongly influenced 
by Kontsevich's idea of {\em motivic integration} \cite{K} and a recent paper
of Denef and Loeser \cite{DL1}. 

\bigskip
\bigskip

\noindent
{\bf Acknowledgements:} It is my  pleasure to thank Professors
Yujiro  Kawamata, Maxim  Kontsevich, Shigefumi Mori, Miles Reid, Yuji 
Shimizu, Joseph  Steenbrink  and  Duco van Straten for useful  
discussions concerning this topic   
and the Taniguchi Foundation for its financial support.

\newpage

\section{Basic definitions and notations}

Let $X$ be a normal irreducible algebraic variety of dimension $d$ over 
$\C$, 
$Z_{d-1}(X)$ the 
group of Weil divisors on $X$,  ${\rm Div }(X) \subset Z_{d-1}(X) $ 
the subgroup of Cartier divisors on 
$X$, $Z_{d-1}(X) \otimes \Q$ 
the group of Weil divisors on $X$ with coefficients in $\Q$, 
$K_X \in Z_{d-1}(X)$ a  
canonical divisor of  $X$. 
We denote by $b_i(X)$ the rank of the $i$-th cohomology 
group with compact supports $H^i_c(X, \Q)$. By the {\bf usual Euler number} 
of $X$ we always mean 
$$ e(X) = \sum_{i \geq 0} (-1)^i b_i(X). $$
The groups   $H^i_c(X, \C)$, $( 0 \leq i \leq 2d)$ 
carry a natural mixed Hodge structure \cite{P-D}. 
We denote by $h^{p,q}\left(H^i_c(X, \C)\right)$ the dimension 
of the $(p,q)$-type Hodge component in $H^i_c(X, \C)$
and define the {\bf $E$-polynomial} by the formula 
\[ E(X; u,v) := \sum_{p,q} e^{p,q}(X) u^p v^q,  \]
where 
\[ e^{p,q}(X): = 
\sum_{i \geq 0} (-1)^i h^{p,q}\left(H^i_c(X, \C)\right). \]
In particular, one has $e(X) = E(X;1,1)$.

\begin{rem} 
{\rm For our purpose, it will be very important that $E$-polynomials 
have properties which are very similar to ones of the  usual Euler 
characteristic: 

(i) if  $X = X_1 \cup \cdots \cup X_k$ is a disjoint union of 
Zariski locally closed strata $X_1, \ldots, X_k$,  then 
\[ E(X; u,v) = \sum_{i=1}^k E(X_i; u,v); \]  

(ii) if $X = X_1 \times X_2$ is a product of two algebraic varieties 
$X_1$ and $X_2$, then 
\[ E(X; u,v) = E(X_1; u,v) \cdot E(X_2; u,v); \]

(iii) if $X$ admits a Zariski locally
trivial fibering over $Z$  such that each fiber of 
$f\,: \, X \to Z$ is isomorphic to affine space ${\C}^n$, then 
\[ E(X; u,v) =  E({\C}^n; u,v) \cdot E(Z; u,v)= (uv)^n E(Z; u,v). \]
}
\label{e-poly} 
\end{rem}

Recall the 
following fundamental definition from the Minimal Model Program \cite{KMM}:

\begin{dfn}
{\rm A variety  $X$ is said to have at worst  
{\bf log-terminal} (resp. {\bf canonical }) singularities, if the 
following conditions are satisfied:   

(i) $K_X$  is an element of ${\rm Div}(X) \otimes \Q$ 
(i.e. $X$ is $\Q$-Gorenstein). 

(ii)  for a  resolution of singularities   
$\rho \,: \, Y  \rightarrow  X$
such that the exceptional locus of $\rho$ is a divisor $D$ 
whose irreducible 
components $D_1, \ldots, D_r$ are  smooth divisors 
with only normal crossings, we have  
\[ K_{Y} = \rho^* K_X + \sum_{i=1}^r a_i D_i  \]
with  $a_i > -1 $ (resp. $a_i \geq  0)$ for all $i$, where $D_i$ runs 
over all irreducible components of $D$. }
\label{log-t} 
\end{dfn}

\section{Stringy Hodge numbers}

\begin{dfn}
{\rm Let $X$ be a normal irreducible  algebraic variety with at worst 
log-terminal  
singularities, $\rho\, : \, Y \rightarrow 
X$ a resolution of singularities  as in \ref{log-t}, 
$D_1, \ldots, D_r$ irreducible components of the exceptional locus, 
and $I: = \{1, \ldots, r\}$. For  any subset 
$J \subset I$ we 
set    
\[ D_J := \left\{ \begin{array}{ll}   
\bigcap_{ j \in J} D_j & \mbox{\rm if $J \neq \emptyset$} 
\\
Y  & \mbox{\rm 
if $J =  \emptyset$} \end{array} \right.  
\;\;\;\; \;\;\;\; \,\mbox{\rm and}   \;\;\;\; \;\;\;\;
D_J^{\circ} := D_J \setminus \bigcup_{ j \in\, I \setminus J} D_j. \]
We define an algebraic  function  
$E_{\rm st}(X; u,v)$ in two 
variables $u$ and $v$
as follows:  
\[ E_{\rm st}(X; u,v) := \sum_{J \subset I} E(D_J^{\circ}; u,v) 
\prod_{j \in J} 
\frac{uv-1}{(uv)^{a_j +1} -1}  \]
(it is  assumed   $\prod_{j \in J}$ to be  $1$, if $J = \emptyset$).
We call  $E_{\rm st}(X; u,v)$ {\bf stringy $E$-function of} $X$. }
\label{def-main}
\end{dfn}

\begin{rem}  

{\rm Since all $a_i \in \Q$, but not necessary $a_i \in \Z$,  
$E_{\rm st}(X; u,v)$ {\bf is not a rational} 
function in  $u,v$ in general. On the other hand, if all singularities 
of $X$ are Gorenstein, then the 
numbers $a_1, \ldots, a_r$ are nonnegative 
integers, hence,   $E_{\rm st}(X; u,v)$ is a {\bf rational 
function}. Moreover, it is easy to see that 
 $E_{\rm st}(X; u,v)$ is an element 
of $\Z[[u,v]] \cap \Q(u,v)$. 
} 
\end{rem}

\begin{dfn}
{\rm In the above situation, we call the rational number  
\[ e_{\rm st}(X) : = \lim_{u,v \to 1} E_{\rm st} (X; u,v) = 
\sum_{J \subset I} e(D_J^{\circ}) \prod_{j \in J} 
\frac{1}{a_j +1} \]
the {\bf stringy Euler number} of $X$.} 
\end{dfn}

For  $E_{\rm st}(X; u,v)$ and $e_{\rm st}(X)$ to be well-defined
we need   the following main  technical statement: 

\begin{theo} 
The function $E_{\rm st}(X; u,v)$ doesn't depend on the choice of 
a resolution $\rho\,: \, Y \rightarrow X$, if all irreducible 
components of the exceptional divisor $D$ are normal crossing divisors.
Moreover, it is  enough to demand the normal crossing condition only for 
those irreducible components $D_i$ of $D$ for which  $a_i \neq 0$ 
$($see the formula in  {\rm \ref{log-t}}$)$. 
\label{key}
\end{theo} 

\begin{rem}
{\rm A complete proof of this statement will be postponed to 
Appendix. Its  main idea, due to Kontsevich,  
is  based on an interpretation 
of the formula for  $E_{\rm st}(X; u,v)$ as a some sort of 
a ``motivic non-Archimedian integral'' over the space of arcs 
$J_{\infty}(Y)$ which   imitates properties of 
$p$-adic integrals  over   $\Z_p$-points of algebraic schemes over $\Z_p$ 
(see \cite{D1}). Some  details of this  idea can be found  in papers of 
Denef and Loeser \cite{DL1,DL2}
(see also \cite{D2,DL}). } 
\end{rem}

\begin{coro} 
If $X$ is smooth, then $ E_{\rm st}(X; u,v) = E(X; u,v)$. 
\label{smooth}
\end{coro}

Our  next statement shows that although  
$E_{\rm st} (X; u,v)$ is   not necessary a polynomial,  
or even a rational function, it satisfies the  Poincar\'e duality
as well as usual $E$-polynomials of smooth projective 
varieties.  

\begin{theo}  {\sc (Poincar\'e duality)} 
Let $X$ be a projective  
${\bf Q}$-Gorenstein algebraic variety of dimension $d$ with at worst 
log-terminal singularities. Then 
$E_{\rm st} (X; u,v)$ has   the  
following  properties:

{\rm (i) } 
  \[ E_{\rm st} (X; u,v) = (uv)^d E_{\rm st} (X; u^{-1},v^{-1});  \]

{\rm (ii)} 
\[  E_{\rm st} (X; 0, 0) = 1. \]
\label{p-d} 
\end{theo}

\noindent
{\em Proof.} (i) Let $\rho\, : \, Y \rightarrow X$ 
be a resolution with the normal crossing exceptional divisors 
$D_1, \ldots, D_r$. Using the stratification 
\[ D_J = \bigcup_{J' \subset J} D_{J'}^{\circ}, \]
we obtain
\[ E (D_J; u,v) =  \sum_{J' \subset J} E ( D_{J'}^{\circ}; u,v). \]
Applying the last equality to  each closed stratum $D_{J'}$, we come 
to the following formula 
\[  E (D_J^{\circ}; u,v) =  \sum_{J' \subset J} (-1)^{|J| -|J'|} 
 E(D_{J'}; u,v). \]
Now we can rewrite the formula for $E(X; u,v)$ as follows
\[ E_{\rm st} (X; u,v) = \sum_{J \subset I} \left(\sum_{J' \subset J}
 (-1)^{|J| -|J'|} 
 E(D_{J'}; u,v) \right)\prod_{j \in J} 
\frac{uv-1}{(uv)^{a_j +1} -1} = \]
\[  = \sum_{J \subset I} E(D_J; u,v) 
\prod_{j \in J} \left( \frac{uv -1}{(uv)^{a_j +1} -1} -1 
\right) 
 =\sum_{J \subset I} E(D_J; u,v) 
\prod_{j \in J} \left( \frac{uv -(uv)^{a_j + 1}}{(uv)^{a_j +1} -1}  
\right) . \] 
It remains to observe   that each term of the last sum satisfies 
the required duality, because  
\[ (uv)^{|J|} \prod_{j \in J}
 \left( \frac{(uv)^{-1} -(uv)^{-a_j - 1}}{(uv)^{-a_j -1} -1} \right) 
= \prod_{j \in J} \left( \frac{uv -(uv)^{a_j + 1}}{(uv)^{a_j +1} -1}  
\right) \]
and 
\[  (uv)^{d - |J|)}E(D_J; u^{-1},v^{-1}) =  E(D_J; u,v) \]
(the last equation  follows from  the Poincar{\'e} duality for smooth 
projective  subvarieties 
$D_J \subset Y$ of dimension $d-|J|$).  

(ii) By substitution  $u=v=0$ in the expression for 
 $E_{\rm st}(X; u,v)$ obtained in the proof of (i), we get  
\[  E_{\rm st} (X; u,v) = E(Y; 0,0). \]
It remains to use the equality  
$E(Y; 0,0) = 1$, which follows from the Poincar\'e duality for 
$Y$.  

\hfill $\Box$

\begin{dfn}
{\rm Let $X$ be a projective  algebraic variety with at worst Gorenstein 
canonical singularities. Assume that $E_{\rm st}(X; u,v) $  is a polynomial: 
$E_{\rm st}(X; u,v) = \sum_{p,q} a_{p,q} u^p v^q. $
Then we define the {\bf stringy
Hodge numbers} of $X$ to be 
\[ h^{p,q}_{\rm st} (X): = (-1)^{p+q} a_{p,q}. \]
}
\end{dfn}

\begin{rem} 
{\rm It follows immediately from \ref{p-d} that if $E_{\rm st}(X; u,v)$ 
is a polynomial, then 
the degree of $E_{\rm st}(X; u,v)$ is equal to $2d$. Moreover, in this 
case one has 
$h_{\rm st}^{0,0}(X) = h_{\rm st}^{d,d}(X) =1$ and 
$h^{p,q}_{\rm st}(X)$ are integers
satisfying   $h_{\rm st}^{p,q}(X)  = 
h_{\rm st}^{q,p}(X)$ $\forall p,q$.}   
\end{rem} 

\begin{conj} 
 Let $X$ be a projective  algebraic variety with at worst Gorenstein 
canonical singularities. Assume that $E_{\rm st}(X; u,v) $  is a polynomial. 
Then all stringy Hodge numbers   $h^{p,q}_{\rm st} (X)$ are nonnegative. 
\end{conj}

Now we consider a geometric  approach to a computation of 
stringy Hodge numbers. 

\begin{dfn} 
{\rm Let $X$ be a projective algebraic variety with at worst canonical
Gorenstein singularities. A birational projective morphism 
$\rho\;: \; Y \rightarrow X$ is called a {\bf crepant desingularization of} 
$X$ if $Y$ is smooth and $\rho^*K_X = K_Y$. } 
\end{dfn}

\begin{theo}
Assume that a ${\Q}$-Gorenstein algebraic variety 
with at worst log-terminal singularities 
$X$ admits a projective birational 
morphism  
$\rho\, : \, Y \rightarrow X$ such that 
$\rho^* K_X = K_{Y}$. Then 
$ E_{\rm st} (X; u,v)  =  E_{\rm st} (Y; u,v)$.  
In particular, if $\rho$ is a crepant desingularization $($i.e. if 
$Y$ is smooth$)$, then 
$E_{\rm st} (X; u,v)  =  E(Y; u,v)$
and therefore 
$h^{p,q}_{\rm st}(X) = h^{p,q}(Y) \;\; \forall p, q.$ 
\label{d-crep}
\end{theo}

\noindent
{\em Proof.} Let $\alpha\, : \, Z \rightarrow Y$ 
be a resolution of singularities of $Y$  such that the exceptional 
locus of the birational morphism $\alpha$ is a divisor $D$ with 
normal crossing irreducible components $D_1, \ldots, D_r$. 
We have  
\[  K_Z = \alpha^* K_Y + \sum_{i =1}^r a_i D_i. \] 
The composition  $\alpha \circ \rho\, : \, Z \rightarrow X$
 is a resolution of singularities of $X$ and we have 
\[  K_Z = \alpha^* (\rho^*K_X) + \sum_{i =1}^r a_i D_i, \] 
since $\rho^* K_X = K_{Y}$. 
Denote by $D'$ the 
exceptional locus of  $\alpha \circ \rho$. Then $Supp\, D' \supset 
Supp\, D$, i.e., if one writes 
\[  K_Z = \alpha^* (\rho^*K_X) + \sum_{i =1}^{s} a_i' D_i', \]
where $D_1', \ldots, D_s'$ are irreducible components of 
$D'$, then $a_j' =0$ for all $j \in \{1, \ldots, s \}$ such that 
$D_j' \subset Supp\, D' \setminus Supp\, D$. By \ref{key}, 
we obtain the same formulas for the $E$-functions of $Y$ and $X$ in terms 
of the desingularizations $\rho \circ \alpha$ and $\alpha$, i.e., 
$E_{\rm st} (X; u,v)  =  E_{\rm st} (Y; u,v)$. The last statement 
of this theorem follows immediatelly from \ref{smooth}. 
\hfill $\Box$

\begin{rem}  
{\rm 
Let $X$ be an $n$-dimensional  Calabi-Yau variety in the sense of \ref{d-cy}.
In \cite{B1} we have shown that   
\[ \frac{d^2}{d u^2} E_{\rm st}(X; u,1)|_{u =1} = \frac{3n^2 - 5n}{12} 
e_{\rm st}(X). \]

This  identity of the above  type for smooth 
varieties $X$ has appeared in the papers of T. Eguchi et. al 
\cite{EHX,EMX} in the connection 
with the Virasoro algebra. 
If $X$ is smooth, then the above identity  follows 
from the Hirzebruch-Riemann-Roch formula (see \cite{LB1,LW}). By 
\ref{d-crep},  this immediately implies the identity for any  $X$ which admits 
a crepant desingularization. 
If $X$ is a $K3$-surface, then  the  identity is equivalent 
to the equality $e(X) =24$. We remark that for smooth 
Calabi-Yau $4$-folds $X$ with $h^{1,0}(X) = h^{2,0}(X) = 
h^{3,0}(X) =0$ the identity is equivalent to the equality 
\[  e(X) = 6( 8 - h^{1,1}(X) + h^{2,1}(X) - h^{3,1}(X)) \]
(see also \cite{W2}). 
} 
\end{rem}

\section{Stringy  $E$-function of toric varieties}

Let us consider the case  when $X$ is a normal  $d$-dimensional 
$\Q$-Gorenstein  
toric variety associated with a rational polyhedral 
fan $\Sigma \subset N_{\R} = N \otimes {\R} $, where  
${N}$ is a free abelian group of rank $d$. Denote by $X_{\sigma}$ 
the torus orbit in $X$ corresponding to a cone $\sigma$ 
($codim_X X_{\sigma} = dim\, \sigma$). Then the property 
of $X$ to be  $\Q$-Gorenstein is equivalent to existence 
of a continious function 
$\varphi_K\, : \, N_{\R} \rightarrow {\R}_{\geq 0}$ satisfying
the conditions 

(i) $\varphi_K (e) = 1$ , if $e$ is a primitive integral generator of 
a $1$-dimensional cone $\sigma \in \Sigma$ 

(ii) $\varphi_K$ is linear on each cone $\sigma \in \Sigma$. 
\bigskip

The following statement is well-known in toric geometry (see e.g. 
\cite{KMM}, Prop. 5-2-2): 

\begin{prop} 
Let $\rho\: : \; X' \to X$ be a toric desingularization of $X$, which 
is defined by a subdivision $\Sigma'$ of the fan $\Sigma$. Then 
the irreducible components $D_1, \ldots, D_r$ of the exceptional divisor $D$ 
of the birational morphism $\rho$ have only normal crossings and they 
one-to-one correspond to primitive integral generators $e_1', \ldots, e_r'$ 
of those $1$-dimensional cones  $\sigma' \in \Sigma'$ which do not belong 
to $\Sigma$. Moreover,  
in the formula 
\[ K_{X'} = \rho^* K_X + \sum_{i=1}^r a_i D_i,  \]
one has $a_i = \varphi_K(e_i') -1$ $\forall i \in \{1, \ldots, r\}$. 
\end{prop} 

\begin{coro} 
Any  normal  $\Q$-Gorenstein  toric variety has at worst log-terminal 
singularities.  
\end{coro} 
      
Let  $\sigma^{\circ}$ be the relative interior of $\sigma$ 
(we put $\sigma^{\circ} = 0$,
 if $\sigma = 0$).  We give the following  explicit formula for the 
function $E_{\rm st}(X; u, v)$:  

\begin{theo} 
\[ E_{\rm st}(X; u,v) = (uv -1)^d \sum_{ \sigma \in \Sigma} 
\sum_{n \in \sigma^{\circ} \cap N} (uv)^{-\varphi_K(n)}.  \] 
\label{st-tor}
\end{theo} 

\noindent 
{\em Proof.} Let  $\rho\: : \; X' \to X$ be a toric desingularization of $X$
defined by a regular simplicial subdivision  
$\Sigma'$ of the fan $\Sigma$,  
 $D_1, \ldots, D_r$  the irreducible components of the exceptional divisor $D$ 
of $\rho$ corresponding 
to primitive integral generators  $e_1', \ldots, e_r' \in N$, and 
$I = \{1, \ldots, r\}$.
First of all we remark that 
the condition  $\bigcap_{j \in J} D_j \neq \emptyset$  
for some nonempty $J \subset I$  is equivalent 
to the fact that 
\[ \sigma_J := 
\{ \sum_{j \in J} \lambda_j e_j' \; | \; \lambda_j \in \R_{\geq 0} \}  \]
is a cone of $\Sigma'$. Let  $\Sigma'(J)$ be the star of the 
cone $\sigma_J \subset \Sigma'$, i.e.,  $\Sigma'(J)$ consists of the 
cones $\sigma' \in \Sigma'$ such that $\sigma' \supset \sigma_J$. 
We denote by  $\Sigma_{\circ}'(J)$ the  subfan of  $\Sigma'(J)$ consisting 
of those cones $\sigma' \in  \Sigma'(J)$ which do not contain any 
$e_i'$ where $i \not\in J$. Then the fan $\Sigma'(J)$ defines the 
toric subvariety $D_J \subset X'$ and the fan  $\Sigma_{\circ}'(J)$ 
defines the open subset $D_J^{\circ} \subset D_J$. The canonical 
stratification by torus orbits 
\[ X' = \bigcup_{ \sigma' \in \Sigma'} X_{\sigma'}' \]
induces the following stratifications 
\[ D_J^{\circ} =  \bigcup_{ \sigma' \in \Sigma_{\circ}'(J)} X_{\sigma'}',  
\;\;\;\;\; \emptyset \neq J \subset I. \]
So we have 
\begin{equation}
 E(D^{\circ}_J; u,v) = \sum_{ \sigma' \in \Sigma_{\circ}'(J)} 
(uv -1)^{d - dim\, \sigma'} \;\;\;\;\; \emptyset \neq J \subset I.
\label{n-emp}  
\end{equation} 
Let $\Sigma'(\emptyset)$ be the subfan of $\Sigma'$  
 consisting 
of those cones $\sigma' \in  \Sigma'$ which do not contain any element 
of $\{ e_1', \ldots, e_r' \}$.  
Then  $\Sigma'(\emptyset)$ defines the canonical stratification of 
\[ X' \setminus D = \sum_{\sigma' \in \Sigma'(\emptyset)} X_{\sigma'}'. \]
In particular, one has 

\begin{equation}
 E(X' \setminus D; u,v) =  \sum_{ \sigma' \in \Sigma'(\emptyset)} 
(uv -1)^{d - dim\, \sigma'}.
\label{emp}
\end{equation} 
The smoothness of $X'$ implies that 
$\sigma_J \cap N$ is a free semigroup with the basis $\{ e_j' \}_{j \in J}$. 
For this reason  we can express  the power series 
\[  \sum_{ n \in \sigma^{\circ}_J \cap N} (uv)^{- \varphi_K(n)} \]
as a product of $|J|$ geometric series:  
\[ \sum_{ n \in \sigma^{\circ}_J \cap N} 
(uv)^{- \varphi_K(n)} = \prod_{j \in J} \frac{ (uv)^{-\varphi_K(e_j')} } 
{ 1 -  (uv)^{-\varphi_K(e_j')}} = 
\prod_{j \in J} \frac{ (uv)^{-a_j -1} } 
{ 1 -  (uv)^{-a_j -1}} = \prod_{j \in J} \frac{1}{(uv)^{a_j +1} -1}.  \]
Hence, 
\begin{equation}
\prod_{j \in J} \frac{ uv -1 }{ (uv)^{a_j + 1}  -1}  = 
(uv -1)^{|J|} \sum_{ n \in \sigma^{\circ}_J \cap N} (uv)^{- \varphi_K(n)}.
\label{factor}  
\end{equation} 
Combining (\ref{n-emp}), (\ref{emp}) and (\ref{factor}), we come to the 
formula 
\begin{equation}
  E_{\rm st}(X; u,v) =  E(X \setminus D; u,v) + 
\sum_{\emptyset \neq 
J \subset I} E(D_J^{\circ}; u,v)(uv -1)^{|J|} \left( 
\sum_{ n \in \sigma^{\circ}_J \cap N} (uv)^{- \varphi_K(n)} \right) = 
\label{formu}
\end{equation}
\[ =  \sum_{ \sigma' \in \Sigma'(\emptyset)} 
(uv -1)^{d - dim\, \sigma'} +  \sum_{\emptyset \neq 
J \subset I} \left(  \sum_{ \sigma' \in \Sigma_{\circ}'(J)} 
(uv -1)^{d + |J|- dim\, \sigma'} \right) \left( 
\sum_{ n \in \sigma^{\circ}_J \cap N} (uv)^{- \varphi_K(n)} \right).  \]
Note that  
\[ \Sigma' = \Sigma(\emptyset) \cup  \bigcup_{\emptyset \neq 
J \subset I} \Sigma_{\circ}'(J). \]
If $ \sigma' = \{ \lambda_{1} e_{i_1}'  + \cdots +  \lambda_{k} e_{i_k}'
\; | \; \lambda_1, \ldots, \lambda_k \in \R_{\geq 0} \}$
is a $k$-dimensional  cone of $\Sigma'(\emptyset)$, then  
the value of $\varphi_K$ is  $1$ on all elements of 
$\{ e_{i_1}', \ldots, e_{i_k}\}$ 
and therefore 
\[ (uv-1)^d  \sum_{n \in (\sigma')^{\circ} \cap N} (uv)^{-\varphi_K(n)}  
= (uv-1)^d \left(\frac{ (uv)^{-1}}{ 1- (uv)^{-1}}\right)^k=  (uv -1)^{d-k}. \]
On the other hand, if 
$\sigma' = \{ \lambda_{1} e_{i_1}'  + \cdots +  \lambda_{k} e_{i_k}'
\; | \; \lambda_1, \ldots, \lambda_k \in \R_{\geq 0} \}$ 
is a $k$-dimensional  cone of $\Sigma'(J)$ $(\emptyset \neq J \subset I)$, 
then 
\[ (uv-1)^d  \sum_{n \in (\sigma')^{\circ} \cap N} (uv)^{-\varphi_K(n)} = 
 (uv-1)^d \left(\frac{ (uv)^{-1}}{ 1- (uv)^{-1}}\right)^{k-|J|} 
\sum_{ n \in \sigma^{\circ}_J \cap N} (uv)^{- \varphi_K(n)} =  \]
\[=  (uv -1)^{d + |J|-k}   
\sum_{ n \in \sigma^{\circ}_J \cap N} (uv)^{- \varphi_K(n)},   \]
since the value of $\varphi_K$ is  $1$ on all elements of 
$\{ e_{i_1}', \ldots, e_{i_k}' \} \setminus \{ e_j' \}_{j \in J}$. 
It remains to apply (\ref{formu}) and use the fact that $\Sigma'$ is a 
subdivision of $\Sigma$, i.e., 
\[ \bigcup_{\sigma \in \Sigma} \sigma^{\circ} \cap N = 
 \bigcup_{\sigma' \in \Sigma'} (\sigma')^{\circ} \cap N. \] 
\hfill $\Box$

\begin{prop} 
If all singularities of a toric variety $X$ are Gorenstein, then 
$E_{\rm st}(X; u,v)$ is a polynomial. 
\end{prop}  

\noindent
{\em Proof.} According to  \cite{R}, $X$ is Gorenstein if and only if 
$\varphi_K(n) \in \Z$ for all $n \in N$. By \ref{st-tor}, 
it is sufficient to prove that 
\[ (uv -1)^d \sum_{n \in \sigma^{\circ} \cap N} (uv)^{-\varphi_K(n)} \]
is a polynomial in  $uv$ for any $\sigma \in \Sigma$. The latter 
follows from the fact that 
\[ (1-t)^{dim\, \sigma} \sum_{n \in \sigma^{\circ} \cap N} 
t^{\varphi_K(n)}  = T(\sigma,t) \]
is a polynomial in $t$ of degree $dim\, \sigma$ \cite{DKh}.
\hfill $\Box$

Recall the following definition introduced by M. Reid \cite{R}(4.2): 

\begin{dfn} 
{\rm Let $X$ be an arbitary $d$-dimensional  
normal $\Q$-Gorenstein toric variety defined by a 
fan $\Sigma$. Denote by $\Sigma^{(d)}$ the set of all $d$-dimensional 
cones in $\Sigma$.  
If $\sigma \in \Sigma^{(d)}$ is an arbitrary  cone, then define 
{\bf shed of $\sigma$}
to be the pyramid
\[ {\rm shed}\, \sigma = \sigma \cap \{ y \in N_{\R}\, : \, 
\varphi_K(y) \leq 1 \}. \]
Furthermore, define 
{\bf shed of $\Sigma$} to be 
\[ {\rm shed}\, \Sigma = \bigcup_{\sigma \in  \Sigma^{(d)}}  
{\rm shed}\, \sigma. \]
}
\end{dfn}

\begin{dfn} 
{\rm Let $\sigma \in \Sigma^{(d)}$ be an arbitrary cone. 
Define $vol(\sigma)$ to be the volume of ${\rm shed}\, \sigma$ 
with respect to the lattice $N \subset N_{\R}$ multiplied by $d!$.  
We set 
\[  vol(\Sigma) := \sum_{\sigma \in  \Sigma^{(d)}}  
vol(\sigma). \]
} 
\end{dfn} 

\begin{rem} 
{\rm It is easy to see that  $vol(\sigma)$ is a positive integer. 
Moreover,  $vol(\sigma)=1$ if and only if 
$\sigma$ is generated by a $\Z$-basis of $N$. } 
\end{rem} 

\begin{dfn} 
{\rm Let $X,X'$ and $X''$ be normal $d$-dimensional $\Q$-Gorenstein toric 
varieties. Assume that we are given two equivariant projective birational 
morphisms $g\,:\, X \to X''$ and $h\, :\, X' \to X''$ such that 
$K_X^{-1}$ is $g$-ample , $K_{X'}$ is $h$-ample, and both $g$ and $h$ 
are isomorphisms in codimension $1$. Then the birational 
rational map $f= h^{-1} \circ g \, : \, X \dasharrow X'$ 
is called a {\bf toric 
flip}. } 
\end{dfn}

The next statement playing important role in termination of toric 
flips is due to M. Reid (see p. 417 in \cite{R}): 

\begin{prop} 
Let $f \, : \, X \dasharrow X'$ be a toric flip of arbitrary 
$\Q$-Gorenstein toric varieties $X$ and $X'$ 
associated with fans $\Sigma$ and $\Sigma'$.
Then
\[  vol(\Sigma) > vol(\Sigma'). \]
\label{v-m}
\end{prop} 

Now we observe the following:

\begin{prop} 
Let $X$ be an arbitrary normal $\Q$-Gorenstein toric variety defined 
by a fan $\Sigma$. Then  the stringy Euler number $e_{\rm st}(X)$ 
equals $vol(\Sigma)$. In particular,  one always has 
$e_{\rm st}(X) \in \Z$. 
\end{prop}

\noindent
{\em Proof.}  By \ref{d-crep}, it is enough to prove the statement 
only for the case of a simplicial fan $\Sigma$ (one can always subdivide 
an arbitrary  fan  $\Sigma$  into a simplicial fan $\Sigma'$ such that 
${\rm shed}\, \Sigma = {\rm shed}\, \Sigma'$).  Let 
$\sigma \in \Sigma$ be a $d$-dimensional 
simplicial cone with primitive lattice generators $e_1, \ldots, e_d$. 
Then $vol(\sigma)$ is equal to the index of the subgroup in $N$ 
generated by  $e_1, \ldots, e_d$. By a direct computation, 
one obtains that 
 \[  \lim_{u,v \to 1} (uv -1)^d \sum_{n \in \sigma^{\circ} \cap N} 
(uv)^{-\varphi_K(n)} =  vol(\sigma). \] 
Now if $\sigma \in \Sigma$ be a $k$-dimensional 
simplicial cone with primitive lattice generators $e_1, \ldots, e_k$ 
$(k <n)$. Then, as in the previous case $k =d$, we compute that 
 \[  \lim_{u,v \to 1} (uv -1)^k \sum_{n \in \sigma^{\circ} \cap N} 
(uv)^{-\varphi_K(n)}=  vol(\sigma) \]
is a positive integer. Therefore, 
 \[  \lim_{u,v \to 1} (uv -1)^d \sum_{n \in \sigma^{\circ} \cap N} 
(uv)^{-\varphi_K(n)} =   vol(\sigma) 
\lim_{u,v \to 1} (uv -1)^{d-k} = 0. \] 
It remains to apply \ref{st-tor}. 
\hfill $\Box$ 

\begin{coro} 
Let $f \, : \, X \dasharrow X'$ be a toric flip of arbitrary $\Q$-Gorenstein  
toric varieties $X$ and $X'$. Then 
\[ e_{\rm st}(X) > e_{\rm st}(X'). \]
\label{e-m}
\end{coro}

\section{Further examples and open questions} 

\begin{exam}
{\rm Let $X_0$ be a $(d-1)$-dimensional 
smooth projective Fano variety with a very ample 
line bundle $L$ such that $L^{\otimes k} \cong K_{X_0}^{-l}$ 
for some positive integers $k,l$. We define 
$X$ to be the cone over $X_0$ in a projective embedding 
defined by $L$, i.e.,  $X$ is obtained by the contraction 
to a singular point $p \in X$ of the section $D$ in  the $\P^1$-bundle 
$Y := \P({\cal O}_{X_0} \oplus L)$ corresponding to the embedding 
${\cal O}_{X_0} \hookrightarrow {\cal O}_{X_0} \oplus L$. 

Let us compute the function $E_{\rm st}(X; u,v)$ 
by considering  the contraction morphism 
$\rho\, : \, Y \rightarrow X$ as a desingularization of $X$. 
By a direct computation, one  obtains  
\[ K_Y = \rho^*K_X + \left( \frac{k}{l} -1 \right) D, \]
i.e., $X$ has at worst log-terminal singularities. 
Since the stratum $X \setminus p = Y \setminus D$ is isomorphic 
to the total space of 
the line bundle $L$ over $X_0$,  we have 
\[  E(X \setminus p ; u,v) = (uv)  E(X_0; u,v). \]
Using the isomorphism $X_0 \cong D$, we obtain     
\[ E_{\rm st}(X; u,v) =   E(X \setminus p ; u,v) + 
\frac{uv -1}{(uv)^{k/l} -1} E(X_0; u,v) = \]
\[ = E(X_0; u,v)\left( uv + \frac{uv -1}{(uv)^{k/l} -1} \right) 
 = \frac{(uv)^{k/l\; +1} -1}{(uv)^{k/l} -1} E(X_0; u,v). \]
Therefore, one has 
\[ e_{\rm st}(X) = \frac{k+l}{k} 
e(X_0). \] 
} 
\label{fano} 
\end{exam} 
\medskip

\begin{exam} 
{\rm It is not true in general that $E_{\rm st}(X; u,v)$ 
is a polynomial even if $X$ has at worst  Gorenstein canonical singularities: 
 Let us take $X_0$ in Example \ref{fano}  to be a smooth quadic
of dimension $d-1 \geq 2$. Then $k = d-1$, $l =1$ and 
\[  E(X_0; u,v) = \left\{ \begin{array}{ll} 
{\displaystyle \frac{(1 + 
(uv)^{\frac{d-1}{2}})((uv)^{\frac{d+1}{2}} -1)}{uv -1}} 
\; & \; 
\mbox{\rm if $d-1$ is even} \\ 
{\displaystyle \frac{(uv)^{d} -1}{uv -1}}\; &\; 
\mbox{\rm if $d-1$ is odd} \end{array} \right. \]
Hence,  for the quadric cone $X$ over $X_0$ we have 
\[  E_{\rm st} (X; u,v) = \left\{ \begin{array}{ll} 
{\displaystyle 
\left( \frac{(uv)^{d} -1}{(uv)^{d-1} -1)} \right)
\left( \frac{(1 + (uv)^{\frac{d-1}{2}})((uv)^{\frac{d+1}{2}} 
-1)}{uv -1} \right) }
\; & \; 
\mbox{\rm if $d-1$ is even,} \\ 
{\displaystyle \left( \frac{(uv)^{d} -1}{(uv)^{d-1} -1)} \right)
\left( \frac{(uv)^{d} -1}{uv -1} \right) } \; &\; 
\mbox{\rm if $d-1$ is odd}. \end{array} \right. \] 
\[  e_{\rm st} (X) = \left\{ \begin{array}{ll} 
{\displaystyle 
\frac{d(d+1)}{d-1} }
\; & \; \mbox{\rm if $d-1$ is even,} \\ 
{\displaystyle  \frac{d^2}{d-1}}  \; &\; 
\mbox{\rm if $d-1$ is odd}. \end{array} \right. \]

The function  $E_{\rm st} (X; u,v)$ is not a polynomial and 
$e_{\rm st}(X) \not\in \Z$  for 
$d \geq 4$. In particular, stringy Hodge numbers of $X$ 
do not exist if $d \geq 4$. If $d =3$, then $E_{\rm st}(X;u,v)$ equals 
$(uv +1)( (uv)^2 + uv +1)$ and  we obtain $h^{1,1}(X) = h^{2,2}(X) =2$, 
$e_{\rm st}(X) =6$.  
}
\end{exam}

\begin{exam} 
{\rm Let $X \subset \C^4$ be a $3$-dimensional hypersurface  defined 
by the equation $x^2 + y^2 + z^2 + t^3 = 0$,  $\rho_0\, : \, 
V_0 \to \C^4$ the blow up of the point 
$q=(0,0,0,0) \in \C^4$, and $X_0 := \rho_0^{-1}(X)$. Then $X_0$ 
is smooth and the exceptional 
locus of the birational morphism $\rho_0\, : \, X_0 \to X$ consists 
of an irreducible divisor $D_0 \subset X_0$ which  is isomorphic 
to a singular quadric $Q \subset \P^3$ defined by the equation 
$z_0^2 - z_1z_2 =0$. Denote by $p_0 \in D_0$ the unique singular point in 
$D_0$. 
Let $\rho_1\, : \, Y  \to X_0$ be the blow up of $p_0$ on $X_0$, 
$D_1$ the birational transform  of $D_0$, and $D_2 = \rho_1^{-1}(p_0)$.  
Then the composition $\rho = \rho_1 \circ \rho_0\; :\; Y \to X$ 
is a resolution of $X$ with normal crossing divisors $D_1, D_2$. 
The unique singularity $q \in X$ is terminal. On the other hand, one has 
\[ K_Y = \rho^*K_X +  1 \cdot D_1 + 2 \cdot D_2, \]
\[ E(D_{\emptyset}; u,v) = (uv)^3 -1, \; E(D_1^{\circ};u,v) = (uv +1)uv, \; 
E(D_2^{\circ};u,v) = (uv)^2, \; E(D_{\{1,2\}}; u,v) = uv +1. \]
Hence, the stringy $E$-function of the Gorenstein variety $X$ 
is the following  
\[ E_{\rm st}(X; u,v)=  (uv)^3 -1 + (uv +1)uv \frac{uv-1}{(uv)^2 -1} 
+ (uv)^2 \frac{uv-1}{ (uv)^3 -1 } + 
(uv +1) \frac{(uv-1)^2}{((uv)^2 -1)((uv)^3 -1 )} = \]
\[ = (uv)^2\frac{(uv)^3 +(uv)^2 + 2uv +1}{(uv)^2 + uv +1}, \]     
i.e., $E_{\rm st}(X; u,v)$ is not a polynomial and $e_{\rm st}(X) = 5/3 
\not\in \Z$.   
}
\end{exam}

\begin{rem} 
{\rm It is known that canonical Gorenstein  singularities of surfaces are 
exactly $ADE$-rational double points  which always 
admit crepant desingularizations. By \ref{d-crep}, 
$E_{\rm st}(X; u,v)$ is a polynomial and 
$e_{\rm  st}(X) \in \Z$ 
for arbitrary algebraic surface $X$ with at worst canonical singularities. 
One can prove that $e_{\rm  st}(X) \in \Z$ also 
for arbitrary algebraic surface $X$ with at worst 
log-terminal singularities. }
\end{rem}

\begin{ques} 
Let $X$ be a geometric quotient of ${\C}^n$ modulo an action of a  
semisimple subgroup $G \subset SL(n, \C)$. Is it   
true that  $E(X; u,v)$ is a polynomial? 
\end{ques} 

\begin{rem} 
{\rm The answer is known to be positive  if  $G$ is commutative ($X$ 
is a Gorenstein toric variety), or if $G$ is finite (see \cite{BD}). }
\end{rem}

\begin{conj} 
Let  $f\; : \; X \dasharrow Y$ be a flip of two projective $3$-dimensional 
varieties with log-terminal singularities. Is it true 
that $e_{\rm st} (X) > e_{\rm st}(Y)$? 
\end{conj} 

\begin{rem} 
{\rm As we observed  in \ref{e-m}, the statement  holds true  
if $X$ and $Y$ are toric varieties and $f$ is a toric flip.}   
\end{rem} 

\begin{conj} 
Let $X$ be a $d$-dimensional algebraic variety with at worst Gorenstein 
canonical singularities. Then the denominator of the rational number 
$e_{\rm st}(X)$ is bounded by a constant $C(d)$ depending only on  
$d$. 
\end{conj}

\begin{rem} 
{\rm In Example \ref{fano}, it  follows from the boundedness of the 
Fano index $k \leq dim\, X_0 + 1 =d$ that $e_{\rm st}(X) \in 
\frac{1}{d!}\Z$. We expect that if $dim\, X = 3$, 
then $e_{\rm st}(X) \in \frac{1}{n} \Z$, where $n \in \{1,2,3,4,6\}$. 
Some evidences for that can be found in \cite{JK}. 
}   
\end{rem}

\section{Appendix: non-Archimedian integrals} 

Let $X$ be a smooth $n$-dimensional complex manifold, $p \in X$ a point, 
and $\Delta_0:= \{z\, : \, |z| < \delta \} \subset \C$ 
a disc with center $0$ of any  small radius $\delta$. Two germs of 
holomorphic mappings $y_1, y_2\, : \, \Delta_0 \to X$ such that 
$y_1(0)=y_2(0)=p$  
 are called {\bf $l$-equivalent}  
if their derivatives in $0$ coincide up to order $l$. The set of 
$l$-equivalent germs is  denoted by 
$J_l(X,p)$ and called the {\bf jet space  of order $l$ at $p$} 
(see  \cite{G-G}, Part A). 
It is well-known  that the union 
$$J_l(X) = \bigcup_{p \in X} J_l(X,p)$$  
is a complex manifold of dimension 
$(l+1)n$, which is a holomorphic affine  bundle over $X$ (but not a vector 
bundle for $l \geq 2$!). The complex manifold 
$J_l(X)$ is  called the {\bf jet space  of order $l$ of  $X$}.   
There are canonical mappings  
$j_l\, : \, 
J_{l+1}(X) \to J_l(X)$ $(l \geq 0)$ whose fibers 
are isomorphic to affine linear spaces $\C^n$. We denote by 
$J_{\infty}(X)$ the projective limit of $J_l(X)$ and by 
$\pi_l$ the canonical projection $J_{\infty}(X) \to J_l(X)$. 
The space  $J_{\infty}(X)$  is known as the {\bf space of arcs of $X$}.    
We denote by $J_{\infty}(X,p)$ the projective limit of $J_l(X,p)$ and call it  
{\bf set of arcs at $p$}.

Let  $R$ be the formal power series ring $\C[[t]]$, i.e., 
the inverse limit of 
finite dimensional $\C$-algebras $R_l: = \C[t]/(t^{l+1})$.
If  $X$ is  $n$-dimensional  smooth quasi-projective algebraic 
variety over $\C$, then the set of points in 
$J_{\infty}(X)$ (resp. $J_l(X)$) coincides with the set of   
$R$-valued (resp. $R_l$-valued) points of $X$. 

From now on  we shall consider only the spaces $J_{\infty}(X)$, 
where $X$ is a smooth quasi-projective algebraic variety. In this 
case, $J_l(X)$ is a smooth   quasi-projective algebraic variety 
for all $l \geq 0$.

Our further terminology is influenced  by 
the theory of Gaussian measures in infinite dimensional linear topological 
spaces (see the book of  Gelfand and Vilenkin \cite{G-V}, Ch. IV). 
 
\begin{dfn} 
{\rm A set $C \subset J_{\infty}(X)$ is called {\bf cylinder set} if 
there exists a positive integer $l$ such that 
$C = \pi^{-1}_l(B_l(C))$ for some constructible subset    
$B_l(C) \subset J_l(X)$ (i.e., for  a finite union 
of Zariski locally closed subsets). Such a constructible subset  $B_l(C)$ 
will be called the 
$l$-{\bf base of} $C$. By definition, the empty 
set $\empty \subset  J_{\infty}(X)$
is a cylinder set  and its $l$-base in $J_l(X)$
is assumed to be empty for all $l \geq 0$.} 
\end{dfn} 

\begin{rem} 
{\rm  Let $C  \subset J_{\infty}(X)$ be a cylinder set with 
an $l$-base $B_l(X)$. 

(i) It is clear that $B_{l+1}(C):= 
j^{-1}_l(B_l(C)) \subset J_{l+1}(X)$ is the $(l+1)$-base of $C$ and 
$B_{l+1}(X)$ is a Zariski locally trivial affine bundle over 
$B_l(C)$, whose fibers are 
isomorphic to the affine space $\C^n$. 

(ii)  Using (i), it is a standard exercise to show that finite unions, 
intersection and complements  of 
cylinder sets are  again cylinder sets. 
\label{cyl}
} 
\end{rem} 

Next statement follows straightforward from the definition of jet spaces 
of order $l$:

\begin{prop} 
Let $p \in X$ be a point of a smooth complex quasi-projective 
algebraic $n$-fold $X$ 
and $z_1, \ldots, z_n$ are local holomorphic coordinates at    
$p$. Then the local coordinate functions 
$z_1, \ldots, z_n$ define  canonical  isomorphisms of quasi-projective 
algebraic varieties 
\[ J_{l}(X,p)  \cong J_{l}({\C}^n, 0) \;\; \forall l \geq 0. \]
In particular, one obtains a canonical isomorphism  of  
algebraic pro-varieties 
\[ J_{\infty}(X,p) \cong  J_{\infty}({\C}^n, 0) \]
which induces a bijection between cylinder sets in $J_{\infty}(X,p)$ 
and cylinder sets in  $J_{\infty}({\C}^n, 0)$, 
where   $J_{\infty}({\C}^n, 0)$ can be identified with the set of all 
$n$-tuples of formal power series $(p_1(t), \ldots, p_n(t)) \in R^n$ having 
the property $(p_1(0), \ldots, p_n(0)) = (0, \ldots, 0)$. 
\label{point} 
\end{prop} 

\begin{coro} 
For any point $c = (c_1, \ldots, c_n) \in \C^n$, there exists a 
canonical isomorphism  of  
algebraic pro-varieties 
\[ J_{\infty}({\C}^n, c) \cong J_{\infty}(\C^n)   \] 
which is defined by the mapping 
\[ \varphi_c \; :\; 
(x_1(t), \ldots, x_n(t)) \mapsto  ((t-c_1)^{-1}x_1(t), \ldots, 
(t-c_n)^{-1} x_n(t)). \]
The mapping $\varphi_c$  induces a bijection between cylinder sets 
in $J_{\infty}(\C^n,c)$ 
and cylinder sets in  $J_{\infty}({\C}^n)$. 
\label{point1} 
\end{coro} 

Recall  the following property of constructible sets (cf. \cite{G-D}, 
Cor. 7.2.6): 

\begin{prop} 
Let $K_1 \supset K_2 \supset \cdots $ be an infinite 
decreasing sequence of constructible subsets of a complex  algebraic 
variety $V$. Assume that $K_i$ is nonempty for all $i \geq 1$. Then 
\[ \bigcap_{i =1}^{\infty} K_i \neq \emptyset. \]
\label{const-b}
\end{prop} 

\noindent
{\em Proof.} Denote by $\overline{K}_i$ the Zariski closure of 
$K_i$. By notherian property, there exists a positive integer  $m$ such 
that $\overline{K}_i = \overline{K}_{i+1}$ for all $i \geq m$. Let 
$Z$ be an irreducible component of $\overline{K}_m$. Then 
for any $i \geq m$ there exists a nonempty 
Zariski open subset $U_i \subset Z$ such that $U_i$ is contained in 
$K_i$. By theorem of Baire,  $\bigcup_{i =1}^{\infty} U_i \neq \emptyset$. 
Therefore,  $\bigcup_{i =1}^{\infty} K_i \neq \emptyset$.

\hfill $\Box$ 

The following property of cylinder sets will be important: 

\begin{theo} 
Assume that a cylinder set $C \subset J_{\infty}(X)$ is contained 
in a countable union $\bigcup_{i =1}^{\infty} C_i$ of 
cylinder sets $C_i$. Then there exists a positive integer  $m$ such that 
$C \subset \bigcup_{i =1}^{m} C_i$. 
\label{cover}
\end{theo} 

\noindent
{\em Proof.} Without loss of generality, we can assume that all $C_i$ are 
contained in $C$ (otherwise we could consider  cylinder sets 
$C \cap C_i$ instead of $C_i$). 
Define new  cylinder sets
\[ Z_k := C \setminus \left( C_1 \cup \cdots \cup C_k \right), 
\; \;( k \geq 1). \]
Thus,  we obtain  a desreasing sequence 
\[ Z_1 \supset Z_2 \supset Z_3 \supset \cdots \]
of  cylinder sets whose intersection is empty. 
Assume that none of sets $Z_1, Z_2, \ldots$ is empty.   
Then  there exits 
an increasing sequence of positive integers  $l_1 < l_2 < l_3 \cdots $ such 
that $Z_k = \pi^{-1}_{l_k} (B_{l_k}(Z_k))$ for some nonempty 
constructible sets $B_{l_k}(Z_k) \subset J_{l_k}(X)$  for all $k \geq 1$. 
Then 
$\pi_l(Z_k) \neq \emptyset$ for all $l \geq 0$, $k \geq 1$.
By theorem of Chevalley, $\pi_l(Z_k) \subset J_l(X)$ is constructible 
for all $l \geq 0$, $k \geq 1$. On the other hand, using \ref{const-b}, 
we obtain  
\[ B_l:= 
\bigcap_{k \geq 1} \pi_l(Z_k) \neq \emptyset, \;\; \forall l \geq 0.  \]
Since $B_0 \neq \emptyset$, we can choose a closed point $p \in B_0$ and  set 
\[ Z_k(p) := Z_k \cap J_{\infty}(X, p). \]
Then the cylinder sets $Z_k(p)$ are nonempty for all $k \geq 1$ and 
 form  a desreasing sequence 
\[ Z_1(p) \supset Z_2(p) \supset Z_3(p)  \supset \cdots \]
whose intersection is empty. Now we show that the last property leads to 
a contradiction. 

By \ref{point}, we can restrict ourselves to the case 
$X = \C^n$ and $p = p_0 \in \C^n$.  
Using the isomorphism  
\[ \varphi_{p_0}\; : \; J_{\infty}({\C}^n, p_0) \cong J_{\infty}(\C^n)\;\; 
\mbox{\rm (see \ref{point1})},  \] 
and the induced  bijection between cylinder sets  
in $J_{\infty}(\C^n,p_0)$ 
and $J_{\infty}({\C}^n)$, we can repeat  
the same arguments for  $J_{\infty}({\C}^n)$ and 
choose  a next point $p_1 \in J_1(\C^n, p_0) \cong J_0(\C^n)$ such that 
$p_1 \in B_1$ and $j_0(p_1) = p_0$. By induction, we  obtain an infinite
sequence of points $p_l \in B_l$ such that 
$j_l(p_{l+1}) = p_l$. It remains to show that the projective limit of 
$\{ p_l\}_{l \geq 0}$ is a point $q \in  
J_{\infty}({\C}^n, p_0)$ which is contained in all cylinder sets $Z_k$ 
$(k \geq 1)$. Indeed, we have $p_{l_l} \in B_{l_k} \subset \pi_{l_k}(Z_k) =
B_{l_k}(Z_k)$. Since $Z_k$ is a cylinder set with the $l_k$-base 
$B_{l_k}(Z_k)$ and $\pi_{l_k}(q) = p_{l_k}$, we obtain  $q \in Z_k$.  
Therefore, $q \in \bigcap_{k \geq 1}  Z_k \neq 
\emptyset$. Contradiction.
 
\hfill $\Box$

\begin{dfn} 
{\rm Let $\Z[\tau^{\pm 1}]$ be the Laurent polynomial 
ring in $\tau$ with coefficients in $\Z$ and  
${A}$ the group algebra 
of $(\Q, +)$ with coefficients in  $\Z[\tau^{\pm 1}]$. 
We denote by $\theta^s \in {A}$ the image of 
$s \in \Q$ under the natural  homomorphism 
$(\Q, +) \to  ({A}^*, \cdot)$, where  ${A}^*$ 
is the multiplicative group of invertible elements in ${A}$ 
($\theta \in {A}$ is considered as a transcendental element over 
$\Z[\tau^{\pm 1}]$). We also write  
$${A}:=  \Z[ \tau^{\pm 1}][\theta^{\Q}] $$ 
and identify the ring ${A}$  with the  direct limit of the subrings
${A}_N :=  \Z[ \tau^{\pm 1}][\theta^{\frac{1}{N}\Z}] \subset 
{A}$, where $N$ runs over the set of all natural numbers.
} 
\end{dfn}  

\begin{dfn} 
{\rm We consider a topology on ${A}$ defined by 
the {\bf non-Archimedian norm}
\[ \| \cdot \| \; : \; {A} \to \R_{\geq 0}  \]
which is uniquely characterised by the properties: 

(i)  $\| ab \|  = \| a\| \cdot \| b\|$, 
$\forall a,b \in  {\cal A}$; 

(ii)  $\| a + b \|  =  \max \{ 
 \| a\|, \| b \| \}$, 
$\forall a,b \in  {\cal A}$ if  $\| a\| \neq  \| b \|$; 

(ii)  $\| a \|  = 1$, $\forall a \in \Z[ \tau^{\pm 1}] 
\setminus \{0 \}$; 

(iii) $\| \theta^s \|  = e^{-s}$ if $s \in \Q$. 

\noindent
The {\bf completion} of ${A}$ (resp. of ${A}_N$) with respect 
to the norm $\| \cdot \|$ will be denoted  
by $\widehat{A}$ (resp. by  $\widehat{A}_N$).  We set 
\[ \widehat{A}_{\infty} := \bigcup_{N \in \N} 
\widehat{A}_N \subset \widehat{A}. \]}
\end{dfn}  

\begin{rem}
{\rm  The noetherian ring 
$\widehat{A}_N$ consists of Laurent power series of $\theta^{1/N}$ with 
coefficients in $\Z[ \tau^{\pm 1}]$. The ring  $\widehat{A}$ 
consists consists of formal infinite sums
\[ \sum_{i=1}^{\infty} a_i \theta^{s_i},\;\;a_i \in \Z[ \tau^{\pm 1}], \]
where $s_1 <  s_2 < \cdots $ is a sequence of rational numbers 
having the property 
\[ \lim_{i \to \infty} s_i = \infty. \]} 
\end{rem}

\begin{dfn} 
{\rm Let  $W \subset V$ is a constructible subset in a 
complex quasi-projective algebraic variety $V$. Assume that 
\[ W = W_1 \cup \cdots \cup W_k \]
is a  union of pairwise nonintersecting Zariski locally closed subsets 
$W_1, \ldots, W_k$. Then we define {\bf $E$-polynomial of} $W$ as 
follows:
\[ E(W; u,v) := \sum_{i=1}^k E(W_i; u,v). \]} 
\label{const-e}
\end{dfn} 

\begin{rem} 
{\rm Using \ref{e-poly}(i), it is easy to check that the above definition
does not depend on the choice of the decomposion of $W$ into a finite 
union of pairwise nonintersecting  Zariski
locally closed subsets.} 
\end{rem}

Now we define a {\bf non-Archimedian cylinder set measure} 
on $J_{\infty}(X)$.

\begin{dfn} 
{\rm  $C \subset J_{\infty}(X)$ be a  cylinder set. 
We define the {\bf non-Archimedian volume} $Vol_X(C) \in 
{A}_1$  of $C$ by the 
following formula: 
\[ Vol_X(C) := 
E(B_l(C); \tau \theta^{-1}, \tau^{-1} \theta^{-1}) \theta^{(l+1)2n} \in 
{A}_1,  \]
where  $C = \pi^{-1}_l(B_l(C))$ and 
$E(B_l(C); u,v)$ is the $E$-polynomial of the $l$-base 
$B_l(C) \subset J_l(X)$. If $C= \emptyset$, we set $Vol_X(C) :=0$.   
}
\end{dfn}

\begin{prop} 
The definition of  $Vol_X(C)$ does not depend on the choice of 
an $l$-base $B_l(C)$. 
\end{prop} 

\noindent
{\em Proof.} By \ref{e-poly}(iii) and \ref{cyl}(i),  
$$E(B_{l+1}(C); 
\tau \theta^{-1}, \tau^{-1} \theta^{-1}) = 
E({\C}^n; \tau \theta^{-1}, \tau^{-1} \theta^{-1}) \cdot 
E(B_l(C); \tau \theta^{-1}, \tau^{-1} \theta^{-1}) =  $$
$$
= \frac{1}{(\theta^2)^n} E(B_l(C); \tau \theta^{-1}, \tau^{-1} 
\theta^{-1}).$$
This implies the required independence on $l$ by induction arguments. 

\hfill $\Box$ 

\begin{prop} 
The non-Archimedian cylinder set measure has the following properties: 

{\rm (i)} If $Z_1$  and $Z_2$ are two cylinder sets such that $Z_1 
\subset Z_2$, then 
\[ \| Vol_X(Z_1) \| \leq  \| Vol_X(Z_2) \|. \]

{\rm (ii)} If $Z_1, \ldots, Z_k$ are cylinder sets, then 
\[ \| Vol_X(Z_1 \cup \cdots \cup Z_k) \| =   \max_{i =1}^k 
\| Vol_X(Z_i) \|. \]

\label{prop12}
\end{prop} 

\noindent
{\em Proof.} (i) It follows 
from \ref{const-e} that for any constructible set $W$ 
one has 
\[ \| E(W; \tau \theta^{-1}, \tau^{-1} \theta^{-1}) \| = e^{2dim\,W}. \]
Therefore 
\[ \| Vol_X(Z_i) \| = 
 \| E(B_l(Z_i); \tau \theta^{-1}, \tau^{-1} \theta^{-1}) \| 
e^{-2n(l+1)}  = e^{2dim\,B_l(Z_i) - 2n(l+1)}, \;\; i =1,2. \]
The inclusion  $Z_1 \subset Z_2$ implies that there exists 
a nonnegative integer $l$ such that $B_l(Z_1) \subset  B_l(Z_2)$, i.e., 
$dim\, B_l(Z_1) \leq dim\, B_l(Z_2)$. This implies the 
statement. 

(ii) Using   our arguments in  (i), it remains to remark that  
\[ dim (W_1 \cup \cdots \cup W_k) =   \max_{i =1}^k \, dim\, W_i. \]

\hfill $\Box$ 

\begin{rem} 
{\rm It follows immediately from \ref{e-poly}(i) that 
\[ Vol_X(C)= Vol_X(C_1) + \cdots + Vol_X(C_k) \]
if $C$ is a finite disjoint union of cylinder sets $C_1, \ldots, 
C_k$. This show that $Vol_X(\cdot)$ is a finitely additive measure. 
The next example shows that an extension of    
$Vol_X(\cdot)$ to a countable additive measure with values 
in $\widehat{A}_1$ needs additional accuracy.}
\label{union} 
\end{rem}  

\begin{exam} 
{\rm Let $C \subset R= \C[[t]]$ be the set consisting of 
all power series $\sum_{i \geq 0} a_i t^i$ such that $a_i \neq 0$ 
for all $i \geq 0$.  
For any $k \in \Z_{\geq 0}$, we define $C_k \subset R$ to 
be the set consisting of 
all power series $\sum_{i \geq 0} a_i t^i$ such that $a_i \neq 0$ 
for all $0 \leq i \leq k$. We identify $R$ with $J_{\infty}({\C})$. 
Then every $C_k \subset  J_{\infty}({\C})$ is a cylinder set 
and  $Vol_{{\C}}(C_k)= (\theta^2 -1)^{k+1}$. Moreover, we have 
\[ C_0 \supset C_1 \supset C_2 \supset \cdots\; , \;\; \mbox{\rm and} \;\; 
C = \bigcap_{k\geq 0} C_k. \]
However, the sequence
\[   Vol_{{\C}}(C_0), \; Vol_{{\C}}(C_1), \;  
Vol_{{\C}}(C_2), \; \ldots \]
does not converge in $\widehat{A}_1$.} 
\end{exam}

\begin{dfn} 
{\rm  We say that a subset $C \subset J_{\infty}(X)$ 
is  {\bf measurable} if  for any positive real number $\varepsilon$
there exists   a sequence of cylinder sets  $C_0(\varepsilon), 
C_1(\varepsilon), C_2(\varepsilon), \cdots $ such that 
$$ \left( C \Delta C_0(\varepsilon) \right) 
\subset \bigcup_{i \geq 1} C_i(\varepsilon) $$ 
and $\| Vol_X(C_i({\varepsilon}))\| < \varepsilon$ for all $i \geq 1$.
If $C$ is measurable, then the element 
\[ Vol_X(C) := \lim_{\varepsilon \to 0} C_0(\varepsilon) 
\in \widehat{A}_1 \]
will be called the {\bf non-Archimedian volume} of $C$. 
}
\end{dfn}

\begin{theo} 
In the above definition, the limit $\lim_{\varepsilon \to 0} 
C_0(\varepsilon)$
exists and  does not depend on the choice of the sequences  
$C_0(\varepsilon), C_1(\varepsilon),  C_2(\varepsilon), \cdots $.   
\end{theo} 

\noindent
{\em Proof.} Let $C$ be a measurable set. For  two positive real numbers 
$\varepsilon$,  $\varepsilon'$  we choose two sequences   
\[ C_0(\varepsilon),  C_1(\varepsilon), C_2(\varepsilon), \cdots 
\;\;\; \mbox{\rm and} \;\;\;  C_0'(\varepsilon'),  C_1'(\varepsilon'), 
C_2'(\varepsilon'), \cdots, \]
such  that  
\[  \| Vol_X(C_i({\varepsilon}))\| < \varepsilon, \;\;\;  
 \| Vol_X(C_i'({\varepsilon'}))\| < \varepsilon' 
\;\;\;(\forall i \geq 1),\]
and
\[ \left( C \Delta C_0(\varepsilon) \right) \subset \bigcup_{i \geq 1} 
C_i(\varepsilon), \;\;\; 
\left( C \Delta C_0'(\varepsilon') \right) 
\subset \bigcup_{i \geq 1} C_i'(\varepsilon'). \]
Then we have 
\[ \left( C_0(\varepsilon) \Delta C_0'(\varepsilon') \right) 
\subset \left( \bigcup_{i \geq 1} C_i({\varepsilon}) \right) 
\cup  \left( \bigcup_{i \geq 1} C_i'({\varepsilon'}) \right). \]
Since   $\left( C_0(\varepsilon) \Delta C_0'(\varepsilon') \right)$ 
is a cylinder set, it follows from \ref{cover} that there exist two 
positive integers $L(\varepsilon)$ and  $L(\varepsilon')$ such that 
\[ \left( C_0(\varepsilon) \Delta C_0'(\varepsilon') \right) 
\subset \left( \bigcup_{i = 1}^{L(\varepsilon)} C_i({\varepsilon}) \right) 
\cup  \left( \bigcup_{i = 1}^{L(\varepsilon')} 
C_i'({\varepsilon'}) \right). \]
By \ref{prop12}, we obtain 
\[ \| Vol_X  \left( C_0(\varepsilon) \Delta C_0'(\varepsilon') \right) \| 
\leq \max \{ \varepsilon, \varepsilon' \}. \] 
Using the inclusions 
\[  C_0(\varepsilon) \setminus 
\left( C_0(\varepsilon) \cap C_0'(\varepsilon') \right)  \subset 
\left( C_0(\varepsilon) \Delta C_0'(\varepsilon') \right), \;\; 
 C_0'(\varepsilon') \setminus 
\left( C_0(\varepsilon) \cap C_0'(\varepsilon')  \right)  \subset 
\left( C_0(\varepsilon) \Delta C_0'(\varepsilon')\right) , \]
we obtain the inequalities
\[   \| Vol_X  \left(  C_0(\varepsilon) \setminus 
\left( C_0(\varepsilon) \cap C_0'(\varepsilon') \right) \right) \|, 
\; \;\;  \| Vol_X  \left(  C_0'(\varepsilon') \setminus 
\left( C_0(\varepsilon) \cap C_0'(\varepsilon') \right) \right) \| 
\leq \max \{ \varepsilon, \varepsilon' \}. \] 
Using \ref{union} and disjoint union decompositions 
$C_0(\varepsilon)=  \left(  C_0(\varepsilon) \setminus 
\left( C_0(\varepsilon) \cap C_0'(\varepsilon') \right) \right) 
\cup \left( C_0(\varepsilon) \cap C_0'(\varepsilon') \right)$ and 
$C_0'(\varepsilon')=  \left(  C_0'(\varepsilon') \setminus 
\left( C_0(\varepsilon) \cap C_0'(\varepsilon') \right) \right) 
\cup \left( C_0(\varepsilon) \cap C_0'(\varepsilon') \right)$, 
we conclude that 
\[  \| Vol_X  C_0(\varepsilon) -  Vol_X  C_0'(\varepsilon') \| 
\leq  \max \{ \varepsilon, \varepsilon' \}. \] 
Now standard arguments show that both limits 
\[  \lim_{\varepsilon \to 0} C_0(\varepsilon), \;\;  
\lim_{\varepsilon' \to 0} C_0(\varepsilon')   \]
exist and coincide. 

\hfill $\Box$ 

The proof of the following statement is a standard exercise: 

\begin{prop} 
Measurable sets possess the following properties: 

{\rm (i)}  Finite unions, finite intersections  of 
measurable sets are measurable. 

{\rm (ii)} If $C$ is a disjoint union of nonintersecting measurable 
sets $C_1, \ldots, C_m$, then 
\[ Vol_X(C) = Vol_X(C_1) + \cdots + Vol_X(C_m). \] 

{\rm (iii)} If $C$ is measurable, then the complement 
$\overline{C}:= J_{\infty}(X) \setminus C$ is measurable.

{\rm (iv)}   If $C_1, C_2, \ldots, C_m, \ldots$ is an infinite sequence 
of nonintersecting measurable sets having the property 
\[ \lim_{i \to \infty} \| Vol_X(C_i) \| =0, \]
then 
\[ C = \bigcup_{i =1}^{\infty} C_i \]
is measurable and 
\[ Vol_X(C) = \sum_{i =1}^{\infty} Vol_X(C_i). \]   
\label{bool}
\end{prop}

\begin{dfn} 
{\rm We shall say that a subset $C \subset J_{\infty}(X)$ 
has  {\bf measure zero} 
if for any positive real number $\varepsilon$ there exists a sequence 
of cylinder sets $C_1(\varepsilon), C_2(\varepsilon), \cdots $ such that 
 $C \subset \bigcup_{i \geq 1} C_i(\varepsilon)$ and  
$\| Vol_X(C_i({\varepsilon}))\| < \varepsilon$ for all $i \geq 1$. } 
\end{dfn}

\begin{dfn} 
{\rm Let $Z \subset X$ be a Zariski closed algebraic  subvariety.   
For any point $x \in Z$, we denote by ${\cal O}_{X,x}$ the ring 
of germs of holomorphic functions at  $x$. 
Let $I_{Z,x} \subset {\cal O}_{X,x}$ be the ideal of germs of holomorphic 
functions vanishing on $Z$. We set 
\[ J_{l}(Z,x) := \{ y \in J_{l}(X,x)\, : \, g(y) =0 \;\; 
\forall \; g \in  I_{Z,x} \}, \;\; l \geq 1,  \]
\[ J_{\infty}(Z,x) := \{ y \in J_{\infty}(X,x)\, : \, g(y) =0 \;\; 
\forall \; g \in  I_{Z,x} \} \]
and 
\[   J_{\infty}(Z) := \bigcup_{x \in Z}  J_{\infty}(Z,x),  
\;\; J_{\infty}(X,Z) := J_{\infty}(X) \setminus  J_{\infty}(Z). \]
The space 
$J_{\infty}(Z) \subset J_{\infty}(X)$ will be called {\bf space of arcs 
with values in $Z$}. 
If $W \subset X$ is a Zariski locally closed subset, 
i.e., $W = Z \cap U$ for some closed $Z$ and open $U$, then we set 
\[  J_{\infty}(X,W) :=  J_{\infty}(Z) \cap   
J_{\infty}(X, X \setminus U). \]} 
\end{dfn} 

\begin{prop} 
Let  $W$ be an arbitrary Zariski locally closed subset in 
a smooth irredicible quasi-projective manifold $X$.
Then $J_{\infty}(X,W) \subset J_{\infty}(X)$ is measurable. Moreover, 
one has  
\[ Vol_X(J_{\infty}(X,W)) = \left\{ \begin{array}{ll} 0
\; & \; 
\mbox{\rm if $dim\, W < dim\, X$} \\ 
Vol_X(J_{\infty}(X))\; &\; 
\mbox{\rm if $dim\, W = dim\, X$.} 
\end{array} \right. \]  
\label{subv}
\end{prop} 

\noindent
{\em Proof.} Using \ref{bool},  one  immediatelly obtains 
that it suffices to prove the statement only in the case 
when $X$ is affine and $W$ is smooth  irreducible 
subvariety of codimension $1$. Then 
$B_l(Z):= \pi_l(J_{\infty}(Z)) \subset J_l(X)$   
is a complete intersection  of $(l+1)$ divisors, i.e., 
$dim\,  B_l(Z) = (n-1)(l+1)$. Let 
$C_l$ be the cylinder sets with the $l$-base $B_l(Z)$. 
Then $C_l$ contains $J_{\infty}(W)$ and $\| Vol_X(C_l) \| = e^{-2(l+1)}$. 
Taking $l$ arbitrary large, we obtain that $\| Vol_X(J_{\infty}(W)) \|= 0$. 

\hfill $\Box$

\begin{dfn}
{\rm  A function $F\, : \, J_{\infty}(X)  \to  
\Q \cup \{\infty\}$ is called 
{\bf measurable}, if  
 $F^{-1}(s)$ is measurable 
for all  $s \in \Q \cup  \{\infty\}$.} 
\end{dfn} 

\begin{dfn} 
{\rm A measurable function  $F\, : \,  J_{\infty}(X) 
\to \Q \cup  \{\infty\}$ is called 
{\bf exponentially integrable}  if  $F^{-1}(\infty)$ has measure zero and 
the series 
\[ \sum_{s \in \Q}  \| Vol_X(F^{-1}(s)) \| e^{-2s} \]
is convergent. 
 If $F$ is exponentially integrable, then 
the sum 
\[  \int_{ J_{\infty}(X)}  e^{-F} := 
  \sum_{s \in \Q}  Vol_X(F^{-1}(s)) \theta^{2s} \in \widehat{A} \]
will be called the {\bf exponential 
integral of ${F}$ over $ J_{\infty}(X) $}. }
\end{dfn}

\begin{dfn} 
{\rm Let $D \subset Div(X)$ be a  subvariety of codimension $1$, 
$x \in D$ a point, and $g \in {\cal O}_{X,x}$ the local equation for 
$D$ at $x$. For any $y \in J_{\infty}(X,D)$, we denote by 
$\langle D, y \rangle_x$ the order of $g(y)$ at $t =0$. The number 
$\langle D, y \rangle_x$ will be called the {\bf intersection number} of 
$D$ and $y$ at $x \in X$.  
 
We define the  function 
\[ F_D\; : \; J_{\infty}(X) \to \Z_{\geq 0}\cup \{\infty\} \]
as follows: 
\[ F_D(y) = \left\{ \begin{array}{ll} 0
\; & \; 
\mbox{\rm if $\pi_0(y) = x \not\in D$} \\ 
\langle D, y \rangle_x \; &\; 
\mbox{\rm if $\pi_0(y) \in D$, but $y \not\in J_{\infty}(D)$} \\
\infty \; &\; 
\mbox{\rm if $y \in J_{\infty}(D)$}
\end{array} \right. \] 
If $D = \sum_{i =1}^m a_iD_i \in Div(X) \otimes \Q$ is a 
$\Q$-divisor (i.e., $a_i \in \Q$ and $D_i$ is irreducible subvariety $\forall 
i \in \{1, \ldots,m \}$), we 
set 
\[ F_D := \sum_{i =1}^m  a_i F_{D_i},  \]
where the symbol  $\infty$ is assumed to 
have  the properties 
$\infty \pm  \infty = \infty$,   $s \cdot \infty =  \infty$ 
($\forall s \in \Q$).
  }
\end{dfn}

\begin{prop} 
Let $D \subset Div(X) \otimes \Q$ be an arbitrary $\Q$-Cartier divisor. 
Then the function 
\[ F_D \;: \; J_{\infty}(X) \to \Q \cup \{\infty\} \]
is mesurable.   Moreover, one has  
\[ \int_{ J_{\infty}(X)}  e^{- F_{D_1 + D_2}}  = 
\int_{ J_{\infty}(X)}  e^{-F_{D_1} - F_{D_2}}\;\;  \;\; 
\forall \; D_1, D_2 \in  
Div(X) \otimes \Q \]
provided both integrals exist. 
\label{func}
\end{prop}

\noindent
{\em Proof.} If  $D \subset Div(X)$ be a  subvariety of 
codimension $1$, then there exists a finite  open covering by affine subsets  
$X = \bigcup_{i =1}^k U_i$ such that $D \cap U_i \subset U_i$ 
is a principal divisor $(f_i)$.  
Since each subset $J_{\infty}(U_i) \subset J_{\infty}(X)$ is a cylinder set, 
it suffices to show  that the $F_D^{-1}(l) \cap J_{\infty}(U_i)$ is  
a cylinder set with some $(l+1)$-base for all $l \in \Z_{\geq 0}$ 
(cf. \ref{cyl}(ii)). The latter follows directly from a  local definition  
of $\langle D, \cdot \rangle$. The rest of the proof follows from the equality:
\[ F^{-1}(s) = \sum_{s_1 + s_2 = s} F_1^{-1}(s_1) \cap 
 F_2^{-1}(s_2), \]
if $F = F_1 + F_2$. 

\hfill $\Box$

\begin{theo} 
Let $\rho\, : \, X' \to X$ be a projective morphism of 
smooth complex varieties, $W = \sum_{i =1}^r e_i W_i$ 
the Cartier divisor defined by the equality 
\[ K_{X'} = \rho^* K_X +  \sum_{i =1}^r e_i W_i.\]
Then  a  function  $F\, : \, J_{\infty}(X) \to  \Q \cup \{\infty \}$ 
is exponentially
integrable if an only if $F \circ \rho  + F_{W}$ is exponentially  
integrable. Moreover, if the latter holds, then 
\[ \int_{J_{\infty}(X)} e^{-F} = \int_{J_{\infty}(X')} 
e^{- F  \circ \rho- F_{W}}. \]
\label{morph}
\end{theo} 

\noindent
{\em Proof.} First we remark that if $C \subset  J_{\infty}(X)$
is a cylinder set,  then  
${\rho}^{-1}(C) \subset  J_{\infty}(X')$ is again a cylinder set. Moreover, 
since $\rho\, : \, X' \to X$ is projective and birational, 
the induced mapping of the spaces of arcs establishes a bijection 
\[ J_{\infty}(X') \setminus  J_{\infty}(W_1 \cup \cdots \cup W_r)   
\; \leftrightarrow \;   
J_{\infty}(X) \setminus J_{\infty}(\rho(W_1  \cup \cdots \cup W_r)). \]
For any integer $k \in \Z \cup \{\infty\}$, we denote by $U_k(X',W)$ the 
cylinder  set consisting of all arcs $y' \in  J_{\infty}(X')$ 
such that $F_{W}(y') = k$. Then  
for any cylinder $C \subset  J_{\infty}(X)$ we have a disjoint union 
decomposition: 
\[ {\rho}^{-1}(C)= \bigcup_{k \in \Z \cup \{\infty\}} 
  {\rho}^{-1}(C) \cap  U_k(X',W). \]
In particular, if we choose $s  \in \Q \cup \{\infty\}$ and set  
\[ C(k, s) : =  {\rho}^{-1}(F^{-1}(s)) \cap  U_k(X',W), \]
then we obtain   
\[ (F \circ \rho)^{-1}(s) =  \bigcup_{k \in \Z \cup \{\infty\}} 
C(k, s)\;\; \mbox{\rm and} \;\;
F^{-1}(s) =  \bigcup_{k \in \Z \cup \{\infty\}} \rho(C(k, s)). \]  
The  key observation is  the fact that  $k= F_W(y')$ 
is the order of vanishing of 
the Jacobian of $\rho$ at any arc  
$y' \in  U_k(X',W)$. Using this fact, we obtain 
\[ Vol_X( \rho(C(k, s))) = \theta^{2k} Vol_{X'}(C(k, s)). \]
Therefore, 
\[ \sum_{s \in \Q \cup \{ \infty \}}  \| Vol_X(F^{-1}(s)) \| e^{-2s} = 
 \sum_{s \in \Q \cup \{ \infty \}}  \sum_{k \in \Z \cup \{\infty\} } 
 \| Vol_X( \rho(C(k, s))) \|  e^{-2s} =  \] 
\[   =  \sum_{s \in \Q \cup \{ \infty \} }  \sum_{k \in \Z \cup \{\infty\} } 
 \| Vol_{X'}( C(k, s) ) \|  e^{-2s +2k} =  
\sum_{s \in \Q \cup \{ \infty \}}  \| Vol_{X'}((\rho \circ F + F_W)^{-1}(s)) 
\| e^{-2s}. \]    
This implies that the sum   
$$\sum_{s \in \Q \cup \{ \infty \}}  \| Vol_X(F^{-1}(s)) \| e^{-2s}$$ 
converges
if and only if the sum 
$$
\sum_{s \in \Q \cup \{ \infty \}}  \| Vol_{X'}((\rho \circ F + F_W)^{-1}(s)) 
\| e^{-2s}$$ converges. In the latter case we obtain the equality 
\[  \sum_{s \in \Q } Vol_X(F^{-1}(s)) \theta^{-2s} =
\sum_{s \in \Q }  Vol_{X'}((\rho \circ F + F_W)^{-1}(s)) 
 \theta^{-2s}. \]

\hfill $\Box$

\begin{theo} 
Let  $D:= a_1D_1 + \cdots + a_r D_r \in Div(X) \otimes \Q$ be 
a $\Q$-divisor. Assume that the subset of all irreducible components 
$D_i$ with  $a_i \neq 0$ form a system of normal crossing divisors. 
Then $F_D$ is exponentially integrable on $J_{\infty}(X)$ if and only 
$a_i > -1$ for all $i \in \{1, \ldots, r \}$. If the 
latter holds, then 
\[  \int_{J_{\infty}(X)} e^{-F_D} = \sum_{J \subset I}  
E(D^{\circ}_J; \tau \theta^{-1},\tau^{-1} \theta^{-1} ) 
(\theta^{-2} -1)^{|J|}\prod_{j \in J} \frac{1}{1 - \theta^{2(1 + a_j)}}
   \]
\label{int-for}
\end{theo} 

\noindent
{\em Proof.} Denote by $J_{\infty}(D) \subset J_{\infty}(X)$ 
the subset of all arcs 
contained in the union $D_1 \cup \cdots \cup D_r$. Since 
$Vol_X(J_{\infty}(D)) =0$, it is sufficient to investigate  
the exponential integral of $F_D$ 
only over $J_{\infty}(X,D):= J_{\infty}(X) \setminus J_{\infty}(D)$. 
For any $(m_1, \ldots, m_r) \in \Z^r_{\geq 0}$, we denote 
\[ U_{m_1, \ldots, m_r}(X,D):= \{ y \in J_{\infty}(X,D)\; : \; 
\langle D_1, y \rangle = m_1, \ldots, \langle D_1, y \rangle = m_r \}. \]
If  $J \subset I:= \{1, \ldots, r\}$, then we set 
\[ U(X, D_J^{\circ} ) := 
\{ y \in J_{\infty}(X,D)\: : \; \pi_0(y) \in D_J^{\circ} \}  \]
and 
\[ M_J:= \{ (m_1, \ldots, m_r) \in  
\Z^r_{\geq 0}\; : \;  m_j > 0 \Leftrightarrow j \in J \}. \]
Thus,  we obtain the following two stratifications of $J_{\infty}(X,D)$: 
\[  J_{\infty}(X,D): = \bigcup_{J \subset I}   U(X, D_J^{\circ} ) \]
and 
\[ J_{\infty}(X,D): =  \bigcup_{(m_1, \ldots, m_r) \in  
\Z^r_{\geq 0} }  U_{m_1, \ldots, m_r}(X,D), \]
where   
$$ U(X, D_J^{\circ}) =  \bigcup_{(m_1, \ldots, m_r) \in  M_J }  
U_{m_1, \ldots, m_r}(X,D).  $$ 
We remark that the value of 
$F_D$ on the cylinder set $U_{m_1, \ldots, m_r}(X,D)$ 
equals $\sum_{j \in J} {2(m_j+1)a_j}$.  
On the other hand, we have 
\[ Vol_X(U_{m_1, \ldots, m_r}(X,D)) = 
E(D^{\circ}_J;\tau \theta^{-1},\tau^{-1} \theta^{-1}) 
\prod_{l =1}^r (\theta^{-2} -1) \theta^{2(m_l+1)} = \]
\[ = 
 E(D^{\circ}_J;\tau \theta^{-1},\tau^{-1} \theta^{-1})
(\theta^{-2} -1)^{|J|} 
\prod_{j \in J} \theta^{2(m_j+1)}, \;\;\mbox{\rm if 
$(m_1, \ldots, m_r) \in  M_J$}. \]
The function $F_D$ is exponentially integrable if and only if the series  
\[  \sum_{J \subset I}  \;\; 
\sum_{(m_1, \ldots, m_r) \in M_J } 
\| Vol_X(U_{m_1, \ldots, m_r}(X,D))\| \prod_{l=1}^r 
e^{-2 (m_j+1)a_j} =  \]
\[ =   \sum_{J \subset I}  \;\; 
\sum_{(m_1, \ldots, m_r) \in M_J } e^{n-|J| +r}
 \prod_{l=1}^r 
e^{-2 (m_j+1)(a_j+1)} \]
is convergent. The latter holds if an only if 
$a_l+1 > 0$ for all $ l \in \{ 1,\ldots, r \}$. 
Using the stratifications $U(X, D_J^{\circ})$ by 
$U_{m_1, \ldots, m_r}(X,D)$, we can rewrite the integral as follows 
\[ \int_{J_{\infty}(X)}  e^{-F_D} = \sum_{J \subset I}  
E(D^{\circ}_J; \tau \theta^{-1},\tau^{-1} \theta^{-1} ) 
(\theta^{-2} -1)^{|J|} \sum_{m \in M_J}  
\prod_{j \in J} \theta^{2(m_j+1)(1 + a_j)}. \]

\hfill $\Box$  
\bigskip
\bigskip

\noindent
{\bf Proof of Theorem \ref{key}}: Let $\rho_1\, : \, Y_1 \to X$ and 
$\rho_2\, : \, Y_2 \to X$ be two resolutions of singularities such that 
\[ K_{Y_1} = \rho_1^* K_X + D_1, \;\; K_{Y_2} = \rho_2^* K_X + D_2 \]
where 
\[ D_1 = \sum_{i=1}^{r_1} a_i' D_i' \; \;\; \mbox{\rm and} \;\;  
D_2 =  \sum_{i=1}^{r_2} a_i'' D_i'', \;\;\;(a_i', a_i'' > -1).   \]
Choosing a resolution of singularities  $\rho_0\, : \, Y_0\to X$ which 
dominates both resolutions $\rho_1$ and $\rho_2$, we obtain 
two morphisms $\alpha_1\, : \, Y_0\to Y_1$ and 
$\alpha_2\, : \, Y_0 \to Y_2$ such that 
$\rho_0 = \rho_1 \circ \alpha_1 =  \rho_2 \circ \alpha_2$.  
We set  $F = F_{ K_{Y_0} -  \rho_0^* K_X}$. Since 
\[ K_{Y_0} -  \rho_0^* K_X =  
( K_{Y_0} - \alpha_i^* K_{Y_i})  + \alpha_i^*D_i, \;\; (i =1,2),  \]
we obtain 
\[  \int_{J_{\infty}(Y_1)} e^{-F_{D_1}} =  
\int_{J_{\infty}(Y_0)} e^{-F} =  \int_{J_{\infty}(Y_2)} e^{-F_{D_2} }\;\;\;
\mbox{\rm (see \ref{morph})}. \]
On the other hand, by \ref{int-for}, we have 
\[ \int_{J_{\infty}(Y_i)} e^{-F_{D_i}} = E_{\rm st}(X; \tau \theta^{-1},  
\tau^{-1} \theta^{-1}), \;\; i \in \{1,2\}.\]
Making the  substitutions 
\[ u = \tau \theta^{-1},  v=  \tau^{-1} \theta^{-1}, \]
we obtain that the definition of the stringy $E$-function 
$E_{\rm st}(X; u,v)$  doesn't depend 
on the choice of resolutions $\rho_1$ and $\rho_2$. 

\hfill $\Box$  

Now we can deduce a Hodge-theoretic version of the results 
in \cite{B}: 

\begin{coro} 
If two smooth $n$-dimensional projective varieties $X_1$ and $X_2$
having trivial canonical classes  are birationally isomorphic, then 
the Hodge numbers of $X_1$ and $X_2$ are the same.
\end{coro} 

\noindent
{\em Proof. } 
Let $\alpha\, : \, X_1 \dasharrow X_2$ be birational map. By Hironaka's 
theorem we can resolve the indeterminancy locus  of $\alpha$ and 
construct some  projective birational morphisms 
\[ \beta_1 \; : \; X \to X_1 \;\; \mbox{\rm and } \;\;  
\beta_2 \; : \; X \to X_2 \]
such that $\alpha = \beta_2 \circ \beta_1^{-1}$. Since the canonical 
classes of both $X_1$ and $X_2$ are  trivial, it follows from \ref{morph} that 
\[ \int_{J_{\infty}(X_1)} e^{0} = \int_{J_{\infty}(X)} 
e^{- F_{W}} =  \int_{J_{\infty}(X_2)} e^{0}, \]
where $W$ is the canonical class of $X$. On the other hand, 
\[ \int_{J_{\infty}(X_i)} e^{0} = Vol_{X_i}(J_{\infty}(X_i)) = 
E(X_i; \tau \theta^{-1}, \tau^{-1} \theta^{-1})\theta^{2n} \;\; ( i =1,2). \]
This implies the equality for the Hodge numbers of $X_1$ and $X_2$. 
 
\hfill $\Box$

\end{document}